\documentstyle[12pt]{article}

\newcommand{\be}{\begin{equation}}
      \newcommand{\ee}{\end{equation}}
      \newcommand{\ba}{\begin{eqnarray}}
       \newcommand{\ea}{\end{eqnarray}}
\newcommand{\ban}{\begin{eqnarray*}}
       \newcommand{\ean}{\end{eqnarray*}}

\newcommand{\pt}{\partial}

 \renewcommand{\o}[2]{\frac{#1}{#2}}
\newcommand{\hf}{\o{1}{2}}

 \newcommand{\qed}{\hspace*{\fill}\rule{3mm}{3mm}\quad}
 \newcommand{\Pf}{\noindent {\em Proof.} }

\newcommand{\Rk}{\noindent {\em Remark.} }

\newcommand{\sect}[1]{\section{#1} \setcounter{equation}{0}}

\newtheorem{theo}{Theorem}[section]

\begin{document}
\newtheorem{lem}[theo]{Lemma}
\newtheorem{prop}[theo]{Proposition}  
\newtheorem{coro}[theo]{Corollary}

\title{A Heat Kernel Lower Bound for Integral Ricci Curvature}
\author{Xianzhe Dai\thanks{Partially supported by NSF Grant \# DMS-9704296} 
  \and Guofang Wei\thanks {Partially supported by NSF Grant \# DMS-9626419.}}
\date{}
\maketitle

\begin{abstract}
In this note we give a heat kernel  lower bound in term of integral Ricci 
curvature, extending Cheeger-Yau's estimate.
\end{abstract}

\newcommand{\inj}{\mbox{inj}}
\newcommand{\vol}{\mbox{vol}}
\newcommand{\diam}{\mbox{diam}}
\newcommand{\Ric}{\mbox{Ric}}
\newcommand{\Iso}{\mbox{Iso}}
\newcommand{\Hess}{\mbox{Hess}} 
\newcommand{\divg}{\mbox{div}}

\sect{Introduction}

Heat kernel is one of the most fundamental quantities in geometry. It can be 
estimated both from above and below in terms of Ricci curvature (see 
\cite{cy,ly, cgt}). The heat kernel upper bound has been extended to integral 
Ricci curvature by Gallot in \cite{g}. Here we extend Cheeger-Yau's lower 
bound \cite{cy} to integral Ricci curvature. 

Our notation for the
integral curvature bounds on a Riemannian manifold $\left( M,g\right) $ is
as follows. For each $x\in M$ let $r\left( x\right) $ denote the smallest
eigenvalue for the Ricci tensor $\mathrm{Ric}:T_{x}M\rightarrow T_{x}M,$ and
for any fixed number $\lambda$
define 
\[ \rho \left( x\right) =\left| \min \left\{ 0,r\left( 
x\right)-(n-1)\lambda\right\} \right|. \]
Then set
\begin{eqnarray*}
k\left( p,\lambda,R\right) &=&\sup_{x\in M}\left( \int_{B\left( x,R\right) } 
\rho^{p}\right) ^{\frac{1}{p}}, \\
\bar{k}\left( p,\lambda, R\right) &=&\sup_{x\in M}\left( 
\frac{1}{\mathrm{vol}B\left(
x,R\right) }\cdot \int_{B\left( x,R\right) } \rho ^{p}\right) ^{\frac{1}{p}}.
\end{eqnarray*}

These curvature quantities evidently measure how much Ricci curvature lies
below $(n-1)\lambda$ in the (normalized) integral sense. And $\bar{k}\left( 
p,\lambda, R\right)=0$ iff $\mathrm{Ric} \geq (n-1)\lambda$. 

Let $ E(x,y,t)$ denote the heat kernel of the Laplace-Beltrami operator on a 
closed manifold $(M, g)$. For any real number $\lambda$ we denote $ E_\lambda 
(\overline{x,y},t)$ the heat kernel on the model space of constant curvature 
$\lambda$. Our main result is

\begin{theo}
Let $n>0$ be an integer, $p> (n+1)/2,$  $\lambda \leq 0$ real numbers  and $D>0$. 
Then there exists an explicitly computable $\epsilon_0=\epsilon(n,p,\lambda, 
D)$ such that for any $(M, g)$ with ${\rm diam} M \leq D$ and $\bar{k}(p, 
\lambda, D) \leq \epsilon_0$,
\[
E(x,y,t) \geq E_\lambda (\overline{x,y},t) - k(p, \lambda, D) C(n, p, \lambda, 
D) (t^{-\o{n+1}{2}}+1),
\]
for any $x, y \in M$ and $t>0$.
\end{theo}

The basic strategy is the same as in Cheeger-Yau, namely, one transplants the 
heat kernel on the model space to $M$ and compare using Duhamel's principle. 
The new difficulty lies in controlling an error term which would be zero in 
the presence of the pointwise Ricci curvature bound. This is overcome by 
employing 1) the mean curvature estimate in \cite{pw}; 2) a comparison of 
volume element (integrated over the directional sphere); 3) Gallot's upper 
bound estimate \cite{g} of the heat kernel, together with a remarkable result 
of Grigor'yan \cite{gr} which furnishes us with a Gaussian upper bound for the 
heat kernel.

\sect{Basic Facts on Heat Kernel}

Here we fix our notation and collect basic facts on the heat kernel which will 
be used in our proof.

As in \cite{cy} we can define the Laplace-Beltrami operator for generalized 
Dirichlet and Neumann boundary conditions on a general Riemannian manifold 
(possibly incomplete) by choosing appropriate domains. The two coincide for a 
complete manifold. The corresponding heat kernel can simply be defined using 
spectral theorem. The heat kernel thus defined is always positive (\cite[Lemma 
1.1]{cy}), which will be essential for our discussion.

The models as used in \cite{cy} need only to have the right mean curvature on 
the distance spheres. Here we restrict our models to the standard models, 
namely, simply connected spaces of constant sectional curvature. The following 
result \cite[Lemma 2.3]{cy} is critical for Cheeger-Yau's theorem as well as 
in here.

\begin{lem} \label{cyt}
Let $E_{\lambda}(r, t)$ denote the heat kernel on the model space of constant 
curvature $\lambda$, where $r=d(x, y)$. Then, for all $r, t >0$, 
\[ \o{\pt}{\pt r} E_{\lambda}(r, t) <0. \]
\end{lem}

As we mentioned before, we also need uniform upper bounds on heat kernel. This 
is established in \cite{g} for integral Ricci curvature. 

\begin{theo} \label{gal}
Given any real number $\lambda \leq 0$, $p>\o{n}{2}$ and $D>0$, there exists an 
explicitly computable $\epsilon_0=\epsilon(n,p,\lambda,D)$ such that for any 
$(M, g)$ with ${\rm diam} M \leq D$ and $\bar{k}(p, \lambda, D) \leq 
\epsilon_0$ one has 
\[ E(x, y, t) \leq C(n,p,\lambda, D)(t^{-p}+1) , \] 
for any $x, y \in M$ and $t>0$.
\end{theo}

However this estimate is not sufficient for our purpose. Fortunately one has 
the following recent amazing result of \cite[Theorem 1.1]{gr}, which translate 
Gallot's estimate into a Gaussian upper bound. We let $f(t), \ g(t)$ denote  
regular functions in sense of \cite{gr} (which includes all piecewise power 
functions with nonnegative exponents).

\begin{theo} \label{grt}
Let $x, y$ be two points on an arbitrary smooth connected Riemannian manifold 
$M$, for which one has
\[ E(x, x, t) \leq \o{1}{f(t)}, \ \ \ \ E(y, y, t) \leq \o{1}{g(t)}, \]
for all $0<t<T\leq \infty$. Then for any $C>4$ and some $\delta=\delta(C)>0$, 
\[ E(x, y, t) \leq \o{4A}{\sqrt{f(\delta t)}\sqrt{g(\delta t)}} e^{-\o{d(x, 
y)^2}{Ct}}, \]
where $A$ is a constant coming from $f$ and $g$.
\end{theo}

\begin{coro} \label{cgr}
With the assumption of Theorem \ref{gal}, we have
\[ E(x, y, t) \leq C(n,p,\lambda, D)(t^{-p}+1) e^{-\o{d(x, y)^2}{5t}}. \]
\end{coro}

The final piece of information we need is a similar Gaussian type estimate on 
the derivative of the heat kernel on the model space.

\begin{prop} \label{edhk}
For the model space we have
\[ |\o{\pt}{\pt r} E_{\lambda}(r, t)| < C(n,\lambda)(t^{-\o{n+1}{2}}+1) 
e^{-\o{d(x, y)^2}{5t}}. \]
\end{prop}

\Pf The key point here is that the space derivative deteriorates the bound
only by a factor of $t^{\hf}$, whereas the time derivative deteriorates the 
bound by a factor of $t$. This follows easily from, say, the gradient estimate
(Harnack inequality) of Li-Yau \cite{ly}, which asserts that a positive 
solution of the heat equation satisfies, for all $\alpha>1$,
\[ \o{|\nabla u|^2}{u^2} - \alpha \o{u_t}{u} \leq \o{n}{\sqrt{2}} 
\o{\alpha^2}{\alpha -1} H + \o{n}{2} \o{\alpha^2}{t} , \]
where $H$ denotes the lower bound on the Ricci curvature.
\qed

\sect{Comparison of Volume Element}

In \cite{pw} a mean curvature comparison estimate is given in terms of 
$k(p,\lambda,R)$ and therefore one obtains the relative volume comparison for 
integral Ricci curvature. Here we need a comparison of integral of the volume 
element just over the directional spheres (instead of the balls).
 
Let $M^n$ be a complete Riemannian manifold and $x\in M$. Around $x$ use 
exponential polar coordinates and write the volume element as 
$d\,\vol = \omega d\theta_{n-1} \wedge dt$, where $d\theta_{n-1}$ is the 
standard volume element on the unit sphere $S^{n-1}(1)$. As $t$ increases 
$\omega$ becomes undefined but we can simply define it to be zero at those 
$t$'s. We have the important equation $\omega' = m\omega$, where the prime 
indicates differentiation along the radial direction and $m$ is the mean 
curvature of the distance spheres around $x$.

In the space form $M_\lambda^n$ of constant sectional curvature $\lambda$, we 
can similarly write the volume element as $d\,\vol = \omega_\lambda 
d\theta_{n-1} \wedge dt$ and $\omega_\lambda' = m_\lambda \omega_\lambda$. 

Define $\psi = \psi(t,\theta) = \max \{0, m(t, \theta)-m_\lambda (t, 
\theta)\}$ and $0$ whenever it becomes undefined. The following mean curvature 
comparison estimate is established in \cite{pw}.

\begin{theo} \label{pwt} For $p>n/2$, $\lambda \leq 0$,
\be
\left(\int_{B\left( x,r\right) }\psi^{2p}d\vol \right)^{\frac 1{2p}}
\leq C\left( n,p\right) \left(k\left( p,\lambda,r \right)\right)^{\frac 12}.
\ee
\end{theo}

With the help of the mean curvature comparison estimate we deduce a comparison 
estimate for the volume element.

\begin{lem} \label{ceve}
There is a constant $C(n,p,\lambda,R)$ such that for any $ p>\frac{n+1}{2}$, $\lambda \leq 0$ 
$r \leq R$, if $k(p,\lambda,R) \leq 1$, then we have
\be
\frac{\int_{S^{n-1}} \omega(r,\theta)d\theta_{n-1}}{\int_{S^{n-1}} 
\omega_\lambda(r,\theta)d\theta_{n-1}} \leq 1+ 
C(n,p,\lambda,R)\left(k(p,\lambda,R)\right)^{\frac 12}
\ee
\end{lem}

\Rk The assumption $k(p,\lambda,R) \leq 1$ is only for the simplicity of the 
statement.

\Pf We will prove a more general relative version.  Define $u(r) 
=\frac{\int_{S^{n-1}} \omega(r,\theta)d\theta_{n-1}}{\int_{S^{n-1}} 
\omega_\lambda(r,\theta)d\theta_{n-1}}$. From the beginning of the proof of 
Lemma 2.1 in \cite{pw}, we have for $0 \leq r_1<r_2\leq R$, 
\[
u(r_2)-u(r_1) \leq
\frac{1}{\vol S^{n-1}} \int_{r_1}^{r_2} \int_{S^{n-1}} \psi 
\frac{\omega}{\omega_\lambda} d\theta_{n-1} \wedge dt.
\]
Using H\"{o}lder's inequality, we have
\ban
\lefteqn{\int_{r_1}^{r_2} \int_{S^{n-1}} \psi \frac{\omega}{\omega_\lambda} 
d\theta_{n-1} \wedge dt } \\
& &\leq \left( \int_0^R \int_{S^{n-1}} \psi^{2p}\omega d\theta_{n-1} \wedge 
dt\right)^{1/2p} \cdot \left( \int_{r_1}^{r_2} \left( 
\omega_\lambda^{-\frac{1}{2p-1}} \int_{S^{n-1}} \frac{\omega}{\omega_\lambda} 
d\theta_{n-1} \right) dt\right)^{1-\frac{1}{2p}} \\
& & \leq C(n,p) \left(k(p,\lambda,R)\right)^{\frac 12} \left( \int_0^R 
\omega_\lambda^{-\frac{1+\alpha}{2p-1}} dt\right)^{\frac{1}{1+\alpha} \cdot 
1-\frac{1}{2p}} \cdot \left( \int_{r_1}^{r_2} \left( \int_{S^{n-1}} 
\frac{\omega}{\omega_\lambda}d\theta_{n-1}\right)^{1+\frac{1}{\alpha}} 
dt\right)^{\frac{\alpha}{\alpha +1} \cdot \left(1-\frac{1}{2p}\right)},
\ean
where $\alpha>0$ is chosen so that $p> \frac{(1+\alpha)n+1}{2}$, therefore 
$\int_0^R \omega_\lambda^{-\frac{1+\alpha}{2p-1}} dt$ is integrable. Thus 
$u(r_2)$ satisfies the integral inequality
\[
u(r_2) - u(r_1) \leq C(n,p,\lambda,R) \left(k(p,\lambda,R)\right)^{\frac 12} 
\left( \int_{r_1}^{r_2} \left(u(t)\right)^{1+\frac{1}{\alpha}} 
dt\right)^{\frac{\alpha}{\alpha +1} \cdot \left(1-\frac{1}{2p}\right)}.
\]
This implies
\[
\left( u(r_2) - u(r_1) \right)_+ \leq C k^{\frac 12} \left( \int_{r_1}^{r_2} 
\left(u(t)\right)^{1+\frac{1}{\alpha}} dt\right)^{\frac{\alpha}{\alpha +1} 
\cdot \left(1-\frac{1}{2p}\right)}.
\]
Let $v=\max \left\{ u - u(r_1),0 \right\}= \left(u - u(r_1)\right)_+$. Then $u 
\leq v+u(r_1)$ and we have
\[
v \leq Ck^{\frac 12} \left( \int_{r_1}^{r_2} \left(v(t) 
+u(r_1)\right)^{1+\frac{1}{\alpha}} dt\right)^{\frac{\alpha}{\alpha +1} \cdot 
\left(1-\frac{1}{2p}\right)}.
\]
Or
\be  \label{ve}
v^{\frac{\alpha+1}{\alpha} \cdot \frac{2p}{2p-1}} \leq \left( Ck^{\frac 
12}\right)^{\frac{\alpha+1}{\alpha} \cdot \frac{2p}{2p-1}} \int_{r_1}^{r_2} 
\left(v(t) +u(r_1)\right)^{\frac{\alpha+1}{\alpha}} dt.
\ee
Write
\[
\left(v(t) +u(r_1)\right)^{\frac{\alpha+1}{\alpha}} = \left[ \left(v(t) 
+u(r_1)\right)^{\frac{\alpha+1}{\alpha} \cdot \frac{2p}{2p-1}} 
\right]^{1-\frac{1}{2p}}.
\]
Now we use the inequality
\[
(a+b)^q \leq 2^{q-1} (a^q+b^q), \ \ a,b \geq 0, \ q \geq 1
\]
to obtain
\[
\left(v(t) +u(r_1)\right)^{\frac{\alpha+1}{\alpha}\frac{2p}{2p-1}} \leq \left[ 
2^{\frac{\alpha+1}{\alpha} \cdot \frac{2p}{2p-1}-1} 
\left(v(t)^{\frac{\alpha+1}{\alpha} \cdot \frac{2p}{2p-1}} 
+u(r_1)^{\frac{\alpha+1}{\alpha} \cdot \frac{2p}{2p-1}} \right) 
\right]^{1-\frac{1}{2p}}.
\]
And letting $w=v^{\frac{\alpha+1}{\alpha} \cdot \frac{2p}{2p-1}}$, (\ref{ve}) 
becomes
\[
w \leq \left( Ck^{\frac 12} \right)^{\frac{\alpha+1}{\alpha} \cdot 
\frac{2p}{2p-1}} \int_{r_1}^{r_2} 2^{\frac{1}{\alpha} +\frac1{2p}} \left( 
w(t)+ u(r_1)^{\frac{\alpha+1}{\alpha} \cdot \frac{2p}{2p-1}} 
\right)^{1-\frac{1}{2p}} dt.
\]
Let $\bar{w}$ be the solution of
\[ \left\{ \begin{array}{ll} \bar{w}' & = 2^{\frac{1}{\alpha} +\frac1{2p}} 
\left( Ck^{\frac 12} \right)^{\frac{\alpha+1}{\alpha} \cdot \frac{2p}{2p-1}} 
\left( \bar{w}(t)+ u(r_1)^{\frac{\alpha+1}{\alpha} \cdot 
\frac{2p}{2p-1}}\right)^{1-\frac{1}{2p}} \\
\bar{w}(r_1) &=0 \end{array} \right.
\]
Then
\[
\bar{w}(r_2) = \left[ \left(u(r_1)^{\frac{\alpha+1}{\alpha} \cdot 
\frac{2p}{2p-1}}\right)^{\frac 1{2p}} + \frac 1{2p}2^{\frac{1}{\alpha} 
+\frac1{2p}} \left( Ck^{\frac 12} \right)^{\frac{\alpha+1}{\alpha} \cdot 
\frac{2p}{2p-1}} (r_2-r_1) \right]^{2p}  -u(r_1)^{\frac{\alpha+1}{\alpha} 
\cdot \frac{2p}{2p-1}}.
\] 
By Gronwall inequality, we have
\[
w \leq \bar{w},\]
which means
\ban
\lefteqn{\left( u(r_2) - u(r_1) \right)_+} \\
&& \leq \left[ \left( u(r_1)^{\frac{\alpha+1}{\alpha} \cdot \frac{1}{2p-1}} + 
\frac1{2p} 2^{\frac{1}{\alpha} +\frac1{2p}} \left( Ck^{\frac 
12}\right)^{\frac{\alpha+1}{\alpha} \cdot \frac{2p}{2p-1}} (r_2-r_1) 
\right)^{2p}-u(r_1)^{\frac{\alpha+1}{\alpha} \cdot \frac{2p}{2p-1}} 
\right]^{\frac{\alpha}{\alpha +1} \cdot \left(1-\frac{1}{2p}\right)}.
\ean
Using the inequality
\[
(x+a)^q-a^q \leq qx(x+a)^{q-1}, \ \ q\geq 1,
\]
we get
\ban
\lefteqn{\left( u(r_2) - u(r_1) \right)_+} \\
&& \leq \left(2^{\frac{1}{\alpha} +\frac1{2p}} 
(r_2-r_1)\right)^{\frac{\alpha}{\alpha +1} \cdot \left(1-\frac{1}{2p}\right)} 
Ck^{\frac 12} \left(u(r_1)^{\frac{\alpha+1}{\alpha} \cdot \frac{1}{2p-1}} 
+\frac1{2p} 2^{\frac{1}{\alpha} +\frac1{2p}} \left( Ck^{\frac 
12}\right)^{\frac{\alpha+1}{\alpha} \cdot \frac{2p}{2p-1}} (r_2-r_1) 
\right)^{q},
\ean
where $q=(2p-1)\frac{\alpha}{\alpha +1} \left(1-\frac{1}{2p}\right)$.

In particular, when $r_1=0$ and $k(p,\lambda,R) \leq 1$, we get
\[
u(r_2) -1 \leq C(n,p,\lambda,R) \left(k(p,\lambda,R)\right)^{\frac 12}.
\]
\qed

\sect{Proof of Theorem}

We follow the same basic strategy as in Cheeger-Yau, starting with the 
Duhamel's
Principle which needs to be justified because of the singularity of the 
distance function at the cut locus.

Using integration by part and heat equation, we have
\ba  \label{EE}
\lefteqn{E(x,y,t) - E_\lambda (\overline{x,y},t) } \nonumber \\&& = 
-\int_0^t\int_{M} \frac{d}{ds}\left( E_\lambda (\overline{x,w},t-s)\right) 
\cdot E(w,y,s) d\,\vol\, ds \nonumber \\ & &
+ \int_0^t\int_{M} E_\lambda (\overline{x,w},t-s)\cdot \frac{d}{ds}\left( 
E(w,y,s)\right) d\,\vol\, ds \nonumber\\
&  & = -\int_0^t\int_{M} \frac{d}{ds}\left( E_\lambda 
(\overline{x,w},t-s)\right) \cdot E(w,y,s) d\,\vol\, ds \nonumber \\ & &
- \int_0^t\int_{M} E_\lambda (\overline{x,w},t-s)\cdot \Delta E(w,y,s) 
d\,\vol\, ds. 
\ea 
Now 
\ban  \frac{d}{ds} E_\lambda & = & -\Delta_\lambda E_\lambda \\
& =& \frac{\partial^2}{\partial r^2} E_\lambda + m_\lambda(r)  
\frac{\partial}{\partial r} E_\lambda \\
& = & -\Delta E_\lambda -\left( m(r,\theta)-m_\lambda (r)\right) 
\frac{\partial}{\partial r} E_\lambda \\
& \leq & -\Delta E_\lambda -\left( m(r,\theta)-m_\lambda (r)\right)_+ 
\frac{\partial}{\partial r} E_\lambda,
\ean
since $\frac{\partial}{\partial r} E_\lambda \leq 0$ by Lemma \ref{cyt}. Hence 
the righthand side of (\ref{EE}) is 
\ban
&\geq & \int_0^t\int_{M} \Delta E_\lambda (\overline{x,w},t-s) \cdot E(w,y,s) 
d\,\vol\, ds \\
&& - \int_0^t\int_{M} E_\lambda (\overline{x,w},t-s)\cdot \Delta E(w,y,s) 
d\,\vol\, ds \\
&& -\int_0^t\int_{M}\left( m(r,\theta)-m_\lambda (r)\right)_+ |\frac{\partial 
E_\lambda }{\partial r} (\overline{x,w},t-s)| \cdot E(w,y,s) d\,\vol\, ds.
\ean
The first two terms combined can be shown to be nonnegative using the same 
argument as in \cite{cy} (using certain convexity property of the distance 
function at the cut locus). The last term is the extra error term, which is
\[
\geq -\int_0^t \left(\int_{M}\left( m(r,\theta)-m_\lambda (r)\right)_+^{q} 
d\,\vol\right)^{1/q}\! \left(\int_{M}\left|\frac{\partial E_\lambda }{\partial 
r} (\overline{x,w},t-s)  E(w,y,s)\right|^{q'} 
d\,\vol\right)^{\frac{1}{q'}}\! ds,\]
for some $q\leq 2p$ to be chosen later. Here $q'=\o{q}{q-1}$.

Now the first factor is controlled by Theorem \ref{pwt} and the volume 
comparison estimate from \cite[Theorem 1.1]{pw}. For the second factor, 
according to Corollary \ref{cgr} and Proposition \ref{edhk}, 
\[ \left|\frac{\partial E_\lambda }{\partial r} (\overline{x,w},t-s)  
E(w,y,s)\right| \leq C [(t-s)^{-\o{n+1}{2}}+1][ s^{-p_1} + 1] 
e^{-\o{d^2(x,w)}{5(t-s)}}e^{- \o{d^2(w,y)}{5s}}. \]
Here $p_1=\o{n}{2} + \alpha$ will be chosen so that $\alpha>0$ is suitably 
small. In order to apply Corollary \ref{cgr} we now need that 
$\bar{k}(\lambda, p_1, D)$ is small than an explicit constant $\epsilon_0$ (as 
determined by Gallot \cite{g}).

We have to deal with the singularity caused by the heat kernel at $t=0$. 
Divide $\int_0^t$ to $\int_0^{t/2} + \int_{t/2}^t$. If $t>1$, then we divide 
further so $\int_0^t=\int_0^{1/2} + \int_{1/2}^{(t-1)/2} +\int_{(t-1)/2}^t$. 
In the latter case the estimate for the middle term is straightforward. For 
$0\leq s\leq t/2$, we have (and we may assume that $s\leq \hf$ by the above 
discussion)
\[ (t-s)^{-\o{n+1}{2}} \leq (t/2)^{-\o{n+1}{2}}, \ \ \ 
e^{-\o{d^2(x,w)}{5(t-s)}} \leq 1, \]
which implies 
\ban
 \int_0^{t/2} \left(\int_M \left|\frac{\partial E_\lambda }{\partial r} 
(\overline{x,w},t-s)  E(w,y,s)\right|^{q'} d\,\vol \right)^{\frac{1}{q'}} ds  
\\
\leq C(t^{-\o{n+1}{2}}+1) \int_0^{t/2} (s^{-p_1}+1)\left(\int_M
e^{- \o{q' d^2(w,y)}{5s}} d\,\vol \right)^{\o{1}{q'}} ds. 
\ean
Now, writing out the integral over the space using the exponential polar 
coordinate around $y$ 
\[ \int_M e^{- \o{q' d^2(w,y)}{5s}} d\,\vol = \int_0^D e^{- \o{q'r^2}{5s}}
(\int_{S^{n-1}}\omega(r, \theta) d\theta) dr. \]
Here it is used that the integral is over the whole manifold.  With the 
curvature assumption on $k(p, \lambda, D)$, we can apply
the comparison estimate for the volume element Lemma~\ref{ceve} and get
\[ \int_M e^{- \o{q' d^2(w,y)}{5s}} d\,\vol \leq C\int_0^D e^{- \o{q'r^2}{5s}}
\omega_{\lambda}(r) dr. \]
We then make a change of coordinate $r_1=\o{r}{\sqrt{s}}$, deducing
\[ \int_M e^{- \o{q' d^2(w,y)}{5s}} d\,\vol \leq Cs^{\o{n}{2}}\int_0^{\infty} 
e^{-q'r_1^2} \o{\omega_{\lambda}(r_1s^{\hf})}{s^{\o{n-1}{2}}} dr_1. \]
Making use of the inequality 
$\o{\omega_{\lambda}(r_1s^{\hf})}{s^{\o{n-1}{2}}} \leq 
r_1^{n-1}e^{(n-1)\sqrt{|\lambda |s}r_1}$ (which can be easily verified) and 
noticing that since $s\leq 1$ this term is dominated by $ e^{-q'r_1^2}$,
we finally arrive at the following estimate for the $\int_0^{t/2}$ part of the 
error term
\[ C(t^{-\o{n+1}{2} -p_1 + \o{n}{2q'} +1}+1), \]
where $p_1$ and $q$ need to be chosen to satisfy the inequality
\[ -p_1 + \o{n}{2q'} +1 > 0. \]

Similarly one has (this time using the exponential polar coordinate around 
$x$)
\[ \int_{t/2}^t \left(\int_M \left|\frac{\partial E_\lambda }{\partial r} 
(\overline{x,w},t-s)  E(w,y,s)\right|^{q'} d\,\vol \right)^{\frac{1}{q'}} ds 
\leq C(t^{-\o{n+1}{2} -p_1 + \o{n}{2q'} +1}+1). \]

Finally, we note that suitable choice for $p_1$ and $q$ can be easily made. 
For example, $q=n+1$ and $p_1=\o{n+1}{2}$ will do.
\qed

Department of Mathematics, University of California, Santa Barbara, CA 93106 

dai@math.ucsb.edu

 wei@math.ucsb.edu

\end{document}